\documentclass[11pt,a4paper,total={6in, 8in}]{article}

\usepackage[colorlinks=true, linkcolor=blue, citecolor=blue, urlcolor=blue]{hyperref}
\usepackage{float}
\usepackage{subcaption}
\usepackage{graphicx}
\usepackage{amsmath}
\usepackage{bm}
\usepackage{rotating}
\newtheorem{thm}{Theorem}

\newtheorem{lm}{Lemma}
\newtheorem{rem}{Remark}

\newtheorem{exmp}{Example}
\usepackage{latexsym}
\usepackage{color}
\usepackage{booktabs}
\usepackage[]{graphicx}
\usepackage{amsfonts}
\usepackage{enumitem}
\usepackage{amsbsy}
\usepackage{physics}

\newcommand{\qed}{\hfill $\Box$}

\begin{document}
	
	\bibliographystyle{plainnat}
	
	\title{Some continuity estimates for ruin probability and other ruin-related quantities}

	\author{by Lazaros Kanellopoulos\footnote{Department of Statistics and Actuarial - Financial Mathematics, University of the Aegean, Samos, Greece. \\
			Email: lkanellopoulos@aegean.gr}}
	
	\setlength{\textwidth}{5.3in} 
	\setlength{\textheight}{8.55in}
	
	\def\baselinestretch{.92}
	\large\normalsize

	\maketitle
	\date
\begin{abstract}
In this paper we investigate continuity properties for ruin probability  in the classical risk model.  Properties of contractive integral operators are used to derive continuity estimates for the deficit at ruin. These results are also applied to obtain desired continuity inequalities in the setting of continuous time surplus process perturbed by diffusion. In this framework, the ruin probability can be expressed as the convolution of a compound geometric distribution $K$ with a diffusion term. A continuity inequality for $K$ is derived and an iterative approximation for this ruin-related quantity is proposed. The results are illustrated by numerical examples. 
\end{abstract}

\vspace{0.5cm}

\vspace{0.5cm}

\section{Introduction}
	\section{Introduction}\label{}

$\quad$ The concept of measuring distances between probability measures is a fundamental one with applications across several areas of mathematics. 
In actuarial sciences, the aim of ruin theory is to model the surplus of an insurance business using a stochastic process and evaluate its ruin probability. This characteristic is  a typical measure for the solvency of a portfolio. An extension of this, which accounts for the severity of ruin, is the distribution of the deficit at ruin (given that ruin occurs). In general, explicit expressions for the ruin probability and/or deficit at ruin are known only in some cases, for instance, when the claim sizes have exponential or phase-type distribution (see \cite{AsmussenAlbrecher2010}, pg 14-15). Therefore, several theoretical approaches have been proposed to approximate, bound, estimate and numerically compute the ruin probability.\par
In insurance mathematics, stochastic models are used to idealize input and output elements to approximate real insurance activities. The problem of stability stated as seeking an appropriate measure of closeness between the ideal and real input element in order to estimate the corresponding deviation in the output. Several theoretical approaches have been developed to analyze optimal choices of metrics between input elements. The problem of stability of the aggregate claim amount over a finite time horizon is formally analyzed in \cite{BeirlantRachev1987}. The estimation of the ruin probability in univariate risk models, using the strong stability method has been investigated in \cite{Kalashnikov2000} and \cite{EnikeevaKalashnikov2001}. Then, the application of this approach has been extended in various directions. For instance,  the stability of the ruin probability in a Markov modulated risk model with L\'evy process with investments was studied in \cite{Rusaityte2001}. A two-dimensional classical risk model was considered in \cite{BenouaretAissani2010} using the strong stability method (see also, \cite{TouaziBenouaret2017}, \cite{BarecheCherfaoui2019} and \cite{HarfoucheBareche2021}). Continuity properties of the surplus process in multidimensional renewal risk models were studied in \cite{GordienkoOrtega2018}. Furthermore, a functional approach that was used to obtain approximations and bounds for some ruin-related quantities was developed in \cite{Politis2006} and \cite{PittsPolitis2008}.  In the discrete-time setup, approximation techniques involving some integral operator were proposed in \cite{GajekRudz2013} and a general formula yielding approximations based on negative binomial mixtures were presented in \cite{SantanRincon2020}.\par
Recently, the Banach contraction principle and fixed point results for contractive operators  have been applied in the theory of risk models. An approximation of the ruin probability under a Markov modulated classical risk model based on the Banach contraction principle was presented in \cite{GajekRudz2018}. It was shown in \cite{Jiang2021} that the ultimate ruin probability can be expressed as the fixed point of a contraction mapping in terms of q-scale functions. Similar contractive approaches have been examined in several recent works (see for instance \cite{GajekRudz2025} and \cite{SanchezBaltazar2020}). Gordienko and V\'azquez-Ortega in \cite{GordienkoOrtega2016} proposed continuity inequalities for ruin probability using properties of contractive mappings. Extending this approach, we derive new continuity bounds for the ruin probability and the deficit at ruin in terms of various choices of probability metrics and we also investigate properties of contractive operators and fixed point results in the context of a continuous time surplus process perturbed by diffusion. To the best of our knowledge, no existing work has employed probability metrics in the setting of a surplus process perturbed by diffusion.\par
The goal in this paper is to propose appropriate probability metrics that allow us to obtain suitable continuity estimates for ruin-related quantities, such as probability of ruin with and without diffusion and the deficit at ruin, in the classical of risk theory. The paper is organized as follows: in the next section we introduce the main quantities of interest in the classical risk model of risk theory and we present the concept of the continuity problem for the ruin probability. In Section 3, we derive a continuity inequality for ruin probabilities and provide an upper bound on the supremum distance between two deficits at ruin including two numerical examples to illustrate the validity of these results. In the Section 4, we derive upper bounds for a ruin-related quantity in the classical risk process perturbed by diffusion and obtain an approximation-using the Banach fixed point theorem (BFPT) to this quantity. The final section summarizes the paper.

\section{Definitions and Preliminaries}
The paper concerns the compound Poisson model with risk process 
\begin{equation}
	\label{risk process}
	U(t)=u+ct-\sum_{n=1}^{N(t)}X_n, \quad t\geq 0,
\end{equation}
where $u\geq 0$ is the initial surplus, $c$ is the premium rate and $N(t)$ denotes the number of claims up to time $t$. The $X_n$'s represent the claim sizes and they are assumed to be independent identically distributed positive random variables (r.v.'s) with a common distribution function  $F$. Also, $N(t)$ is a Poisson process with rate $\lambda$, independent of $X_n$. We further assume $c=\lambda (1+\theta)EX$, where $\theta>0$ is the relative security loading. Let
\begin{equation}
	\label{ruin probability}
	\psi(u):=Pr(\underset{t\geq 0}{\inf} \; U(t) <0|U(0)=u), \quad u\geq 0,
\end{equation}
be the ruin probability in infinite time. In the classical risk model, it is well-known that $\psi(u)$ satisfies 
\begin{equation}
	\label{tail-G(x)}
	\psi(u)=Pr(L>u)=\sum_{n=1}^{\infty}\frac{\theta}{1+\theta} \left(\frac{1}{1+\theta}\right)^n \overline{F_e^{*n}}(u), \quad u\geq 0,
\end{equation}
where $\overline{F_e^{*n}}(u)=Pr(X_e^{(1)}+X_e^{(2)}+\cdots X_e^{(n)}> u)$ is the tail of the $n$th-fold convolution of $F_e(u)=\int_{0}^{u}\overline{F}(z)dz/\int_{0}^{\infty}\overline{F}(z)dz$ with itself. The variable $L$ here is the maximal aggregate loss in the surplus process. The distribution $F_e$ is known as the equilibrium distribution associated of $F$. \par
In general, the continuity problem is based on the following implication. Let $U(t, \alpha)$ be a risk process governed by a parameter $\alpha=(\lambda, c, F)$ with the ruin probability $\psi_\alpha$. If $\mathcal{A}$ denotes the space of possible values of the parameter then one can view the ruin probability as a mapping $\psi: \mathcal{A} \rightarrow \Psi$, where $\Psi$ is the functional space of all possible functions $\psi_\alpha$. Assume that $\mathcal{A}$ and $\Psi$ are metric spaces with metrics $\delta$ and $\nu$, respectively. In such terms, the problem of interest is reduced to investigation of $(\delta, \nu)$-continuity of this mapping: that is if the implication 
\[
\delta(\alpha, \widetilde{\alpha})\rightarrow 0, \quad \mbox{then} \quad \nu(\psi_\alpha, \psi_{\widetilde{\alpha}})\rightarrow 0,
\]
for $\alpha, \widetilde{\alpha} \in \mathcal{A}$, holds at a fixed point $\alpha$.\par
The metrics $\delta$ and $\nu$ should be computationally convenient and reflect the core of the problem. In the following, we consider various choices for these metrics. Other such metrics that may be suitable to obtain bounds of this form for ruin probabilities and ruin-related quantities can be found in \cite{Kalashnikov1997}.
If we can find an inequality
\begin{equation}
	\label{continuity-inequality}
	\nu (\psi_\alpha, \psi_{\widetilde{\alpha}})\leq w(\delta(\alpha, \widetilde{\alpha})),
\end{equation}
where $w\geq 0, w(s)\overset{s\rightarrow 0}{\rightarrow}0$ and $w(0)=0$, then the inequality $\eqref{continuity-inequality}$ is called a continuity inequality (estimate) and provides the possibility of bounding $\nu(\psi_\alpha, \psi_{\widetilde{\alpha}})$ in terms of the distance $\delta(\alpha, \widetilde{\alpha})$.\par
In applications, the parameters $\lambda$ and $F$ that govern the risk model are usually unknown. Therefore, the intensity of claim arrivals $\lambda$ and the distribution of claim sizes $F$ are approximated by some parameter $\widetilde{\lambda}$ and distribution $\widetilde{F}$, respectively, for which the ruin probability $\widetilde{\psi}$ can be found. This allows local estimations of the ruin probability alterations to be obtained with respect to disturbances of the parameters of $F$ and $\lambda$. Throughout the paper, we make the assumption that $\mu:=\int_{0}^{\infty}xdF(x)<\infty$ and $\widetilde{\mu}:=\int_{0}^{\infty}xd\widetilde{F}(x)<\infty$ hold, as well as the net profit conditions:
\begin{equation}
	\label{net profit}
	\phi=\frac{\lambda \mu}{c}<1, \quad \mbox{and} \quad \widetilde{\phi}=\frac{\widetilde{\lambda} \widetilde{\mu}}{c}<1.
\end{equation}
\section{Continuity estimate for ruin probability and for deficit at ruin in the classical risk model}
\subsection{Ruin probability $\psi(u)$}
\label{subsection-Ruin}
In this section, we will study continuity conditions of the ruin probability between surplus processes as presented in $\eqref{risk process}$. Firstly, we introduce the background needed to use properties of contractive operators in order to bound the ruin probability in terms of an appropriate distance. \par
Let $\mathcal{D}_\infty$ be the set of functions $h:[0, \infty]\rightarrow \mathbb{R}$ which are right-continuous on $[0, \infty)$ and have left-hand limits on $[0, \infty]$ (c\'adl\'ag functions). This space endowed with the supremum norm, $||h||_\infty= \sup_{t \in \mathbb{R}}|h(t)|$, is a nonseperable Banach space (see e.g. \cite{Pollard1984}). For $f\in \mathcal{D}_\infty$, we define a new function $h_\gamma: [0, \infty]\rightarrow [0, \infty]$ such that 
\[
h_\gamma(x)=(1+x)^\gamma |h(x)|, \quad x \geq 0.
\]
Then, for $\gamma\geq 0$, the space of functions $\mathcal{D}_\gamma=\{h : [0, \infty) \rightarrow \mathbb{R}: h_\gamma \in \mathcal{D}_\infty\}$ with the norm on $\mathcal{D}_\gamma$ defined by $||h||_\gamma=||h_\gamma||_\infty$ is again a nonseperable Banach space (see e.g. \cite{GrubelPitts1993}). For any $\gamma\geq 0$, let $\mathcal{D}_\gamma$ endowed with the metric $\nu_{\gamma}$ where, for $x, y \in \mathcal{D}_\gamma$,
\begin{equation}
	\label{weighted-function-metric}
	\nu_\gamma(x, y):=\int_{0}^{\infty}(1+t)^\gamma |x(t)-y(t)|dt, \quad \gamma \geq 0.
\end{equation}
Then $(\mathcal{D}_\gamma, \nu_\gamma)$ is a complete metric space.
If the ruin probability has a power law decay (e.g. large claims with Pareto tails) then it is reasonable to consider a metric $\nu_{\gamma}(\psi_\alpha, \psi_{\widetilde{\alpha}})$. \par
Suppose now that the claim size distribution has finite moments of order $\gamma+1$ for $\gamma\geq 0$, so that $F\in \mathcal{D}_\gamma$. Then the following theorem holds.

\begin{thm}
	\label{LKdistance}
	Let $X$ and $\widetilde{X}$ be two r.v.'s representing the individual claim amounts in two surplus processes as in $\eqref{risk process}$ with distribution functions $F$ and $\widetilde{F}$, respectively. We assume that the intensities of claim arrivals are $\lambda$ and $\lambda_X$ and that $EX^{(\gamma+1)}<+\infty$ and $E\widetilde{X}^{(\gamma+1)}<+\infty$ for some $\gamma\geq 0$.   We denote by $\psi(u)$ and $\widetilde{\psi}(u)$ the ruin probabilities associated with the surplus processes $\{U(t):t\geq 0 \}$ and $\{\widetilde{U}(t):t\geq 0\}$, respectively. Then,
	\begin{align}
		\label{metric-theorem}
		\int_{0}^{\infty}(1+t)^\gamma \left|\psi(t)-\widetilde{\psi}(t)\right|dt &\leq \frac{c}{c-\lambda M^X_\gamma}   \left( \frac{\nu_{\gamma+1}(F, \widetilde{F})}{\gamma+1}+\nu_{\gamma}(F, \widetilde{F}) M^L_\gamma \right. \nonumber \\
		& \left. +\frac{|\lambda-\widetilde{\lambda}|}{c} M^{\widetilde{X}}_{\gamma+1} (1+M^L_\gamma) \right), 
	\end{align}
	where,
	\[
	M^X_\gamma=\int_{0}^{\infty}(1+t)^\gamma\overline{F}_X(t)dt \quad \mbox{and provided that} \quad \lambda M^X_\gamma/c < 1.
	\]
\end{thm}

\textbf{Proof:}\\
	The following integral equation for ruin probabilities is commonly known 
	\begin{equation}
		\label{ruin-1}
		\psi(u)=\frac{\lambda}{c}\left(\int_{u}^{\infty}\overline{F}(t)dt+\int_{0}^{u}\psi(u-t)\overline{F}(t)dt\right).
	\end{equation}
	It is easy to check that for the operators
	\begin{equation}
		\label{operator-1}
		Tx(u)=\frac{\lambda}{c}\left(\int_{u}^{\infty}\overline{F}(t)dt+\int_{0}^{u}x(u-t)\overline{F}(t)dt\right), \quad u\geq 0,
	\end{equation}
	and
	\begin{equation}
		\label{operator-2}
		\widetilde{T}x(u)=\frac{\widetilde{\lambda}}{c}\left(\int_{u}^{\infty}\overline{\widetilde{F}}(t)dt+\int_{0}^{u}x(u-t)\overline{\widetilde{F}}(t)dt\right), \quad u\geq 0,
	\end{equation}
	we have $T\mathcal{D}_\gamma\subset \mathcal{D}_\gamma$ and $\widetilde{T}\mathcal{D}_\gamma \subset \mathcal{D}_\gamma$.\\
	By $\eqref{operator-1}$ and $\eqref{operator-2}$ for every $x,y \in \mathcal{D}_\gamma$,
	\begin{align*}
		\nu_\gamma(Tx, Ty)&=\frac{\lambda}{c}\int_{0}^{\infty}(1+u)^\gamma \left| \int_{0}^{u}x(u-t)\overline{F}(t)dt-\int_{0}^{u}y(u-t)\overline{F}(t)dt\right|du\\
		&\leq \frac{\lambda}{c}\int_{0}^{\infty}\int_{0}^{u}(1+u)^\gamma \overline{F}(t)\left| x(u-t)-y(u-t)\right|dtdu\\
		&= \frac{\lambda}{c}\int_{0}^{\infty}\int_{t}^{\infty}(1+u)^\gamma \overline{F}(t)\left| x(u-t)-y(u-t)\right|dudt,
	\end{align*}
	changing the limits of integration we get the last equality. Hence, setting $u=t+z$ and keeping in mind inequality $(1+z+t)^\gamma\leq (1+z)^\gamma  (1+t)^\gamma$, it follows that
	\begin{align*}
		\nu_\gamma(Tx, Ty)&\leq \frac{\lambda}{c}\int_{0}^{\infty}\int_{t}^{\infty}(1+u)^\gamma \overline{F}(t)\left| x(u-t)-y(u-t)\right|dudt\\
		&=\frac{\lambda}{c}\int_{0}^{\infty}\int_{0}^{\infty}(1+z+t)^\gamma \overline{F}(t)\left| x(z)-y(z)\right|dzdt\\
		&\leq\frac{\lambda}{c}\int_{0}^{\infty}(1+t)^\gamma \overline{F}(t) \int_{0}^{\infty}(1+z)^\gamma \left| x(z)-y(z)\right|dzdt\\
		&=\nu_\gamma(x,y) \frac{\lambda}{c} \int_{0}^{\infty}(1+t)^\gamma \overline{F}(t)dt=\nu_\gamma(x,y) \frac{\lambda}{c} M^X_\gamma.
	\end{align*}
	Therefore, these operators are contractive on $\mathcal{D}_\gamma$ with modules $\lambda \mu/c$ and $\widetilde{\lambda} \widetilde{\mu}/c$, respectively, since net conditions hold (see equation $\eqref{net profit}$)
	Now, integrating by parts, we get
	\begin{align*}
		M^X_\gamma=\int_{0}^{\infty}(1+t)^\gamma \overline{F}(t)dt
		&=\left[\overline{F}(t)  \frac{(1+t)^{\gamma+1}}{\gamma+1}\right]_0^\infty +\int_{0}^{\infty}f(t) \frac{(1+t)^{\gamma+1}}{\gamma+1} dt\\
		&=\frac{E(X+1)^{\gamma+1}-1}{\gamma+1}.
	\end{align*}
	
	Also, it is well known that $E(|X|^p)<\infty \Leftrightarrow E(|X-a|^p)<\infty$, $\forall a \in \mathbb{R}$ and $0<p<\infty$. Therefore, $M^X_\gamma$ exists if the moments $EX^{(\gamma+1)}$ of the claim-size distribution are finite for $\gamma\geq 0$.\\
	According to $\eqref{ruin-1}$, $\psi$ and $\widetilde{\psi}$ are the unique fixed points of $T$ and $\widetilde{T}$ that is  $\psi=T\psi$ and $\widetilde{\psi}=\widetilde{T} \widetilde{\psi}$. Now,
	\begin{align*}
		\nu_\gamma(\psi, \widetilde{\psi})&=\nu_\gamma(T\psi, \widetilde{T}\widetilde{\psi})\\
		&\leq \nu_\gamma(T\psi, T\widetilde{\psi})+\nu_\gamma(T\widetilde{\psi}, \widetilde{T}\widetilde{\psi})\\
		&\leq \frac{\lambda}{c}M^X_{\gamma} \nu_\gamma(\psi, \widetilde{\psi})+\nu_\gamma(T\widetilde{\psi}, \widetilde{T}\widetilde{\psi}),
	\end{align*}
	or,
	\begin{equation}
		\label{proof-1}
		\nu_\gamma(\psi, \widetilde{\psi})\leq \frac{c}{c-\lambda M^X_{\gamma}} \nu_\gamma(T\widetilde{\psi}, \widetilde{T} \widetilde{\psi}),
	\end{equation}
	where $M^X_\gamma=\int_{0}^{\infty}(1+t)^\gamma\overline{F}(t)dt=\frac{1}{\gamma+1}\left(E(X+1)^{\gamma+1}-1\right)$. In view of $\eqref{operator-1}$ and $\eqref{operator-2}$, for each $\psi \in \mathcal{D}_\gamma$, we have
	\begin{align*}
		\nu_\gamma(T\psi, \widetilde{T}\psi)
		&\leq \int_{0}^{\infty}(1+u)^\gamma \left|\frac{\lambda}{c}\left(\int_{u}^{\infty}\overline{F}(t)dt+\int_{0}^{u}\psi(u-t)\overline{F}(t)dt\right) \right.\\
		&\left.-\frac{\lambda}{c}\left(\int_{u}^{\infty}\overline{\widetilde{F}}(t)dt+\int_{0}^{u}\psi(u-t)\overline{\widetilde{F}}(t)dt\right)\right|du\\
		&+\int_{0}^{\infty}(1+u)^\gamma \left|\left(\frac{\lambda-\widetilde{\lambda}}{c}\right)\left(\int_{u}^{\infty}\overline{\widetilde{F}}(t)dt+\int_{0}^{u}\psi(u-t)\overline{\widetilde{F}}(t)dt\right)\right|du\\
		&=I_1+I_2.
	\end{align*}
	For the first term $I_1$ on the last inequality, we have
	\begin{align*}
		I_1&\leq \frac{\lambda}{c}\int_{0}^{\infty}(1+u)^\gamma \left(\int_{u}^{\infty}\left|\overline{F}(t)-\overline{\widetilde{F}}(t)\right|dt+\int_{0}^{u}\psi(u-t)\left|\overline{F}(t)-\overline{\widetilde{F}}(t)\right|dt\right)du\\
		&=\frac{\lambda}{c} \left[\int_{0}^{\infty}\int_{0}^{t}(1+u)^\gamma\left|\overline{F}(t)-\overline{\widetilde{F}}(t)\right|dudt \right.\\ &\left.+\int_{0}^{\infty}\int_{t}^{\infty}(1+u)^\gamma\psi(u-t)\left|\overline{F}(t)-\overline{\widetilde{F}}(t)\right|dudt \right],
	\end{align*}
	where for the first term inside the square bracket it follows that
	\begin{align}
		\label{eq-1}
		\int_{0}^{\infty} \left|\overline{F}(t)-\overline{\widetilde{F}}(t)\right| \int_{0}^{t}(1+u)^\gamma dudt&=\int_{0}^{\infty}\frac{(1+t)^{\gamma+1}-1}{\gamma+1}\left|\overline{F}(t)-\overline{\widetilde{F}}(t)\right|dt \nonumber \\ 
		&\leq \frac{1}{\gamma+1} \nu_{\gamma+1}(F, \widetilde{F}),
	\end{align}
	and
	\begin{align}
		\label{eq-2}
		\int_{0}^{\infty}\int_{t}^{\infty}(1+u)^\gamma\psi(u-t)\left|\overline{F}(t)-\overline{\widetilde{F}}(t)\right|dudt&=\int_{0}^{\infty}\left|\overline{F}(t)-\overline{\widetilde{F}}(t)\right| \nonumber \\
		&\times \int_{0}^{\infty}(1+z+t)^\gamma\psi(z)dz dt  \nonumber \\
		&\leq \int_{0}^{\infty}(1+t)^\gamma \left|\overline{F}(t)-\overline{\widetilde{F}}(t)\right|dt  \nonumber \\
		&\times \int_{0}^{\infty}(1+z)^\gamma \psi(z)dz \nonumber \\ \leq \nu_\gamma(F, \widetilde{F}) M^L_\gamma .
	\end{align}
	Therefore, by $\eqref{eq-1}$ and $\eqref{eq-2}$ it follows that
	\begin{equation}
		\label{proof-2}
		I_1\leq  \frac{1}{\gamma+1} \nu_{\gamma+1}(F, \widetilde{F})+\nu_\gamma(F, \widetilde{F}) M^L_\gamma.
	\end{equation}
	Similarly, for the term $I_2$, we have
	\begin{align*}
		I_2&=\frac{|\lambda-\widetilde{\lambda}|}{c}\left(\int_{0}^{\infty}\overline{\widetilde{F}}(t)\int_{0}^{t}(1+u)^{\gamma}dudt	+\int_{0}^{\infty}\overline{\widetilde{F}}(t)\int_{t}^{\infty}(1+u)^{\gamma}\psi(u-t)dudt\right)\\
		&=\frac{|\lambda-\widetilde{\lambda}|}{c}\left(\int_{0}^{\infty}\overline{\widetilde{F}}(t) \frac{(1+t)^{\gamma+1}-1}{\gamma+1}dt+ \int_{0}^{\infty}\overline{\widetilde{F}}(t)\int_{0}^{\infty}(1+z+t)^\gamma \psi(z)dz dt\right)\\
		&\leq \frac{|\lambda-\widetilde{\lambda}|}{c}\left(\int_{0}^{\infty}\overline{\widetilde{F}}(t) \frac{(1+t)^{\gamma+1}}{\gamma+1}dt+M^L_{\gamma} \int_{0}^{\infty}\overline{\widetilde{F}}(t)(1+t)^\gamma dt\right)\\
		&\leq \frac{|\lambda-\widetilde{\lambda}|}{c}\left(\int_{0}^{\infty}\overline{\widetilde{F}}(t)(1+t)^{\gamma+1}dt+M^L_{\gamma} \int_{0}^{\infty}\overline{\widetilde{F}}(t)(1+t)^{\gamma+1}dt\right)
	\end{align*}
	Hence,
	\begin{equation}
		\label{proof-3}
		I_2\leq \frac{|\lambda-\widetilde{\lambda}|}{c} M^{\widetilde{X}}_{\gamma+1}\left( 1+ M^L_{\gamma}\right).
	\end{equation}
	By $\eqref{proof-1}, \eqref{proof-2}$ and $\eqref{proof-3}$ the result follows.
\hfil \qed \\

\begin{rem}
	For $\gamma=0$, the inequality $\eqref{metric-theorem}$ reduces to the following inequality  
	\[
	\mathbb{K}(\psi, \widetilde{\psi}) \leq \frac{c}{c-\lambda \mu} \left(\nu_1(F, \widetilde{F})+\frac{\mathbb{K}(F, \widetilde{F})EX^2}{2\theta \mu}+\frac{|\lambda-\widetilde{\lambda}|\widetilde{\mu}}{c} \left(1+\frac{EX^2}{2\theta \mu}\right)\right),
	\]
	where $M_\gamma^X=EX=\mu$, $M_\gamma^L=EL=EX^2/2\theta EX$  and the function $\mathbb{K}: \mathcal{F} \times \mathcal{F} \rightarrow [0, \infty]$ with $\mathbb{K}(F,\widetilde{F})=\int_{0}^{\infty}\left|F(t)-\widetilde{F}(t)\right|dt$ is called the \textit{Kantorovich metric} (see e.g. \cite{Kalashnikov1997}) in $\mathcal{F}$ the set of distribution functions $F$ of all positive r.v.'s. 
\end{rem}

\begin{exmp}
	Let the r.v. $X$ follows a mixture of two exponential distributions with tail
	\[
	\overline{F}_{X}(t)=\frac{1}{2} e^{-5t/4}+\frac{1}{2}e^{-5t/6}, \quad t>0,
	\]
	and $\widetilde{X}\sim Exp(1)$ such that $EX=E\widetilde{X}=1$. For convenience, we refer to the continuity bound in $\eqref{metric-theorem}$ as $DK1$. In Table $\ref{LKdistance-table}$ we compare the results of $DK1$ with the real values of distance $\nu_{\gamma}(\overline{F}, \overline{\widetilde{F}})$.
\end{exmp}

\begin{table}[h]
	\caption{Comparison between the bound $DK1$ in \eqref{metric-theorem} and the exact values of $\nu_{\gamma}(F, \widetilde{F})$ for a mixture of two exponential distributions versus an exponential distribution.}
	\label{LKdistance-table}
	\centering
	
	\begin{minipage}{0.48\textwidth}
		\centering
		\textbf{(a) $\gamma=0,\ \lambda=\tfrac{5}{6}$} \\[2mm]
		\begin{tabular}{|c|c|c|c|}
			\hline
			& $c=3$ & $c=5$ & $c=7$ \\
			\hline
			$\lVert \psi-\widetilde{\psi} \rVert_{1,\gamma}$  & 0.0154 & 0.0080 & 0.0060 \\
			\hline
			DK1  & 0.1211 & 0.0999 & 0.0929 \\
			\hline
		\end{tabular}
	\end{minipage}
	\hfill
	\begin{minipage}{0.48\textwidth}
		\centering
		\textbf{(b) $\gamma=0,\ \lambda=\tfrac{10}{11}$} \\[2mm]
		\begin{tabular}{|c|c|c|c|}
			\hline
			& $c=3$ & $c=5$ & $c=7$ \\
			\hline
			$\lVert \psi-\widetilde{\psi} \rVert_{1,\gamma}$  & 0.0174 & 0.0089 & 0.0059 \\
			\hline
			DK1 & 0.1271 & 0.1024 & 0.0944 \\
			\hline
		\end{tabular}
	\end{minipage}
	
	\vspace{4mm}
	
	\begin{minipage}{0.48\textwidth}
		\centering
		\textbf{(c) $\gamma=1,\ \lambda=\tfrac{5}{6}$} \\[2mm]
		\begin{tabular}{|c|c|c|c|}
			\hline
			& $c=3$ & $c=5$ & $c=7$ \\
			\hline
			$\lVert \psi-\widetilde{\psi} \rVert_{1,\gamma}$ & 0.0736 & 0.0353 & 0.0231 \\
			\hline
			DK1 & 0.6340 & 0.3548 & 0.2922 \\
			\hline
		\end{tabular}
	\end{minipage}
	\hfill
	\begin{minipage}{0.48\textwidth}
		\centering
		\textbf{(d) $\gamma=1,\ \lambda=\tfrac{10}{11}$} \\[2mm]
		\begin{tabular}{|c|c|c|c|}
			\hline
			& $c=3$ & $c=5$ & $c=7$ \\
			\hline
			$\lVert \psi-\widetilde{\psi} \rVert_{1,\gamma}$ & 0.0850 & 0.0396 & 0.0257 \\
			\hline
			DK1 & 0.7512 & 0.3795 & 0.3047 \\
			\hline
		\end{tabular}
	\end{minipage}
	
\end{table}

\subsection{The deficit at ruin}
\label{section-deficit}
In this section, we apply the contraction mapping method to obtain a bound for the deficit at ruin in classical risk model.  Contractive operator techniques to derive bounds for the deficit at ruin, though this approach requires the existence of the adjustment coefficient (light-tailed claims) were used in \cite{GajekRudz2025}. Gordienko and V\'azquez-Ortega in \cite{GordienkoOrtega2016} proposed a simple technique for continuity estimation for ruin probability in the compound Poisson risk model. In a similar manner, for $y>0$ fixed, we study the comparison between two defective tails of the deficit at ruin $\overline{G}(u,y)$ and $\overline{\widetilde{G}}(u,y)$, respectively. \par
The deficit at ruin (given that ruin occurs), which was introduced in the paper of \cite{GerberGoovaertsKaas1987}, represents the probability that, starting with a surplus $u$, the deficit at ruin does not exceed $y$, i.e. $G(u, y)=Pr\left(|U_T|\leq y, T<\infty|U(0)=u\right)$. If we set $\overline{G}(u,y)=\psi(u)-G(u,y)$ for $u\geq 0$ and $y\geq 0$, then $\overline{G}(u, y)$ satisfies the defective renewal equation (see \cite{Willmot2002}) 
\begin{equation}
	\label{deficit-1}
	\overline{G}(u,y)=\frac{\lambda}{c}\left[\int_{0}^{u}\overline{G}(u-t,y)\overline{F}(t)dt+\int_{u+y}^{\infty}\overline{F}(t)dt\right],
\end{equation}
which is the associated tail of the defective distribution of $G(u,y)$ and we note that it satisfies $\lim_{y \rightarrow \infty}=\psi(u)<1$.\\
Similarly to Section $\ref{subsection-Ruin}$, we consider a surplus process as $\eqref{risk process}$,  $\widetilde{U}(t)$ , but one that is governed by the parameters $\widetilde{\lambda}$ and $\widetilde{F}$ with associated defective tail of the deficit at ruin which satisfies the following integral equation 
\begin{equation}
	\label{deficit-2}
	\overline{\widetilde{G}}(u,y)=\frac{\widetilde{\lambda}}{c}\left[\int_{0}^{u}\overline{\widetilde{G}}(u-t,y)\overline{\widetilde{F}}(t)dt+\int_{u+y}^{\infty}\overline{\widetilde{F}}(t)dt\right]
\end{equation}
In the following theorem we obtain a bound for the uniform metric between two deficits at ruin. Subsequently, numerical examples are provided for various values of $u$ and $y$.
\begin{thm}
	\label{theorem-deficit}
	With the above notation, it holds that
	\begin{equation}
		\label{KanDist-}
		\sup_{u \geq 0}\left|\overline{G}(u,y)-\overline{\widetilde{G}}(u,y)\right| \leq \frac{1}{c-\lambda \mu} \left[\lambda \mathcal{Q}_y(F, \widetilde{F})+|\lambda-\widetilde{\lambda}|\widetilde{\mu}\right],
	\end{equation}
	where $\mathcal{Q}_y(F, \widetilde{F})=\int_{y}^{\infty}|\overline{F}(t)-\overline{\widetilde{F}}(t)|dt$.
\end{thm}

\textbf{Proof:}\\
	Let $\mathcal{X}$ be the space of all functions $x:[0, \infty)\times[0, \infty) \rightarrow [0, 1]$ endowed with the \textit{uniform metric} $\nu_d(x,y):=\sup_{t\geq 0}|x(t)-y(t)|$. Also, let the variable $z$ as a constant so that $\overline{G}(x,z)=\overline{\Gamma}_z(x)$ is a function of $x$. Then $(\mathcal{X}, v)$ is easy to see that it is a complete metric space.\\
	We have $T^z\mathcal{X}\subset \mathcal{X}$ and $\widetilde{T^z}\mathcal{X}\subset \mathcal{X}$ for the following operators
	\begin{equation}
		\label{operator-3}
		T^zx(u)=\frac{\lambda}{c}\left[\int_{z+u}^{\infty}\overline{F}(t)dt+\int_{0}^{u}x(u-t)\overline{F}(t)dt\right], \quad u\geq 0,
	\end{equation}
	and
	\begin{equation}
		\label{operator-4}
		\widetilde{T^z}x(u)=\frac{\widetilde{\lambda}}{c}\left[\int_{z+u}^{\infty}\overline{\widetilde{F}}(t)dt+\int_{0}^{u}x(u-t)\overline{\widetilde{F}}(t)dt\right], \quad u\geq 0.
	\end{equation}
	Therefore, by $\eqref{operator-3}$, for every $x, y, \in \mathcal{X}$, we have
	\begin{align*}
		\nu_d(T^zx, T^zy)&=\frac{\lambda}{c}\underset{u \geq 0}{sup}\left|\int_{0}^{u}x(u-t)\overline{F}(t)dt-\int_{0}^{u}y(u-t)\overline{F}(t)dt\right|\\
		&\leq \frac{\lambda}{c}\underset{u \geq 0}{sup}\int_{0}^{u}\overline{F}(t)\underset{s\in [0, u]}{sup}|x(s)-y(s)|dt\\
		&\leq \frac{\lambda}{c} v_d(x,y) \int_{0}^{\infty}\overline{F}(t)dt=\frac{\lambda \mu}{c} \nu_d(x,y)
	\end{align*}
	The inequality 
	\begin{equation}
		\label{v(Tx,Ty)}
		\nu_d(\widetilde{T^z}x,\widetilde{T^z}y)\leq \frac{\widetilde{\lambda} \widetilde{\mu}}{c} \nu_d(x,y)
	\end{equation}
	is similarly verified.\\
	According to $\eqref{deficit-1}$ and $\eqref{deficit-2}$, $\Gamma_z(x)$ and $\widetilde{\Gamma}_z(x)$ are the unique fixed points of $T^z$ and $\widetilde{T^z}$ that is $\Gamma_z=T^z\Gamma_z$ and $\widetilde{\Gamma}_z=\widetilde{T^z}\widetilde{\Gamma}_z$. Thus, it follows that
	\begin{align*}
		\nu_d(\Gamma_z, \widetilde{\Gamma}_z)&=\nu_d(T^z\Gamma_z, \widetilde{T^z}\widetilde{\Gamma}_z)\\
		&\leq \nu_d(T^z\Gamma_z, T^z\widetilde{\Gamma}_z)+\nu_d(T^z\Gamma_z, \widetilde{T^z}\widetilde{\Gamma}_z)\\
		&\leq \frac{\lambda \mu}{c}  \nu_d(\Gamma_z, \widetilde{\Gamma}_z)+\nu_d(T^z\Gamma_z, \widetilde{T^z}\widetilde{\Gamma}_z),
	\end{align*}
	or equivalently,
	\[
	\nu_d(\Gamma_z, \widetilde{\Gamma}_z)\leq \frac{c}{c-\lambda \mu} \nu_d(T^z\widetilde{\Gamma}_z, \widetilde{T^z}\widetilde{\Gamma}_z).
	\]
	In view of $\eqref{operator-1}$ and $\eqref{operator-2}$ for each $x\in \mathcal{X}$ we have
	\begin{align*}
		\nu_d(T^z x, \widetilde{T^z} x)&\leq \sup_{u \geq 0} \left| \int_{u+z}^{\infty}\overline{F}(t)dt+\int_{0}^{u}x(u-t)\overline{F}(t)dt\right.\\
		&\left. -\frac{\lambda}{c}\int_{u+z}^{\infty}\overline{\widetilde{F}}(t)dt+\int_{0}^{u}x(u-t)\overline{\widetilde{F}}(t)dt \right| \\
		&+\sup_{u \geq 0}\left|\frac{\lambda-\widetilde{\lambda}}{c}\left(\int_{u+z}^{\infty}\overline{\widetilde{F}}(t)dt+\int_{0}^{u}x(u-t) \overline{\widetilde{F}}(t)dt\right)\right|\\
		&\leq \frac{\lambda}{c}\int_{z}^{\infty}|\overline{F}(t)-\overline{\widetilde{F}}(t)|dt+\frac{|\lambda-\widetilde{\lambda}|}{c}\int_{0}^{\infty}\overline{\widetilde{F}}(t)dt\\
	\end{align*}
	
	Combining the last inequality with $\eqref{v(Tx,Ty)}$, we obtain
	\[
	\nu_d(\overline{G}(u,y), \overline{\widetilde{G}}(u,y))\leq \frac{1}{c-\lambda \mu} \left[\lambda \mathcal{Q}_y(F, \widetilde{F})+|\lambda-\widetilde{\lambda}|\widetilde{\mu}\right],
	\]
	where $\mathcal{Q}_y(F, \widetilde{F})=\int_{y}^{\infty}|\overline{F}(t)-\overline{\widetilde{F}}(t)|dt$.
	To prove the similar inequality with $\frac{1}{c-\widetilde{\lambda} \widetilde{\mu}}$ instead of $\frac{1}{c-\lambda \mu}$ inequality $\eqref{v(Tx,Ty)}$ is used and the result follows.
\hfil \qed \\

\begin{exmp}
	\label{Example-deficit}
	Let three non-negative r.v.'s $X^{(1)}\sim Erlang(3,3), X^{(2)}\sim Exp(1)$ and $X^{(3)}$ with tail $\overline{F}_{X^{(3)}}(t)=\frac{1}{2} e^{-5t/4}+\frac{1}{2}e^{-5t/6}$. We also assume the defective distribution function of the deficit at ruin $G^{(1)}(u,y), G^{(2)}(u,y)$ and $G^{(3)}(u,y)$ in risk models with individual claim sizes $X^{(1)}, X^{(2)}$ and $X^{(3)}$, respectively. \\
	In Table $\ref{comparison G(u,y)}$ we present the exact values of the distance $\left|\overline{G}^{(i)}(u,y)-\overline{G}^{(j)}(u,y)\right|,$ $ i,j=1,2,3, i\neq j$ versus the bound $DK2$\ derived from  $\eqref{KanDist-}$ for several cases. For each case, we use the same $\lambda$ (see \cite{WillmotLin1998} and \cite{Willmot2002}). 
\end{exmp}

\begin{table}[htbp]
	\caption{Comparison between $\left|\overline{G}^{(i)}(u,y)-\overline{G}^{(j)}(u,y)\right|$, $i,j=1,2,3$, $i\neq j$ and $DK2$.}
	\label{comparison G(u,y)}
	\centering
	
	\begin{minipage}{0.48\textwidth}
		\centering
		\scriptsize
		\textbf{(a) $\theta = 1$, $X^{(1)}$ vs $X^{(2)}$} \\[2mm]
		\begin{tabular}{c|c|c|c}
			\hline
			$y$ & $u$ & $\left|\overline{G}^{(1)}-\overline{G}^{(2)}\right|$ & $DK2$ \\
			\hline
			0.10 & 0.10 &  0.0080 & 0.2900\\
			& 0.25 &  0.0184 & \\
			& 0.50 &  0.0367 & \\
			& 1.00 &  0.0640 & \\
			& 2.00 &  0.0754 & \\
			&      &         & \\
			0.25 & 0.10 & 0.0227 & 0.2726 \\
			& 0.25 & 0.0365 & \\
			& 0.50 & 0.0550 & \\
			& 1.00 & 0.0742 & \\
			& 2.00 & 0.0736 & \\
			&      &        & \\
			0.50 & 0.10 & 0.0495 & 0.2220 \\
			& 0.25 & 0.0625 & \\
			& 0.50 & 0.0758 &  \\
			& 1.00 & 0.0821 & \\
			& 2.00 & 0.0098 & \\
			&      &        & \\
			1.00 & 0.10 & 0.0773 & 0.1547 \\
			& 0.25 & 0.0818 & \\
			& 0.50 & 0.0832 &  \\
			& 1.00 & 0.0748 &  \\
			& 2.00 & 0.0521 & \\
			&      &        & \\
			2.00 & 0.10 & 0.0521 & 0.1905 \\
			& 0.25 & 0.0506 & \\
			& 0.50 & 0.0463 &   \\
			& 1.00 & 0.0372 &  \\
			& 2.00 & 0.0233 & \\
			\hline
		\end{tabular}
	\end{minipage}
	\hfill
	\begin{minipage}{0.48\textwidth}
		\centering
		\scriptsize
		\textbf{(b) $\theta = 4$, $X^{(1)}$ vs $X^{(2)}$} \\[2mm]
		\begin{tabular}{c|c|c|c}
			\hline
			$y$ & $u$ & $\left|\overline{G}^{(1)}-\overline{G}^{(2)}\right|$ & $DK2$ \\
			\hline
			0.10 & 0.10 &  0.0034 & 0.0735 \\
			& 0.25 &  0.0084 & \\
			& 0.50 &  0.0176 & \\
			& 1.00 &  0.0290 & \\
			& 2.00 &  0.0241 & \\
			&      &        & \\
			0.25 & 0.10 & 0.0091 & 0.0681 \\
			& 0.25 & 0.0100 & \\
			& 0.50 & 0.0233 & \\
			& 1.00 & 0.0304 & \\
			& 2.00 & 0.0228 & \\
			&      &        & \\
			0.50 & 0.10 & 0.0194 & 0.0555 \\
			& 0.25 & 0.0243 & \\
			& 0.50 & 0.0294 &  \\
			& 1.00 & 0.0303 & \\
			& 2.00 & 0.0194 & \\
			&      &        & \\
			1.00 & 0.10 & 0.0300 & 0.0387 \\
			& 0.25 & 0.0309 & \\
			& 0.50 & 0.0301 &  \\
			& 1.00 & 0.0247 &  \\
			& 2.00 & 0.0132 & \\
			&      &        & \\
			2.00 & 0.10 & 0.0205 & 0.0476 \\
			& 0.25 & 0.0189 & \\
			& 0.50 & 0.0161 &   \\
			& 1.00 & 0.0114 &  \\
			& 2.00 & 0.0054 & \\
			\hline
		\end{tabular}
	\end{minipage}
	
	\vspace{4mm}
	
	\begin{minipage}{0.48\textwidth}
		\centering
		\scriptsize
		\textbf{(c) $\theta = 4$, $X^{(1)}$ vs $X^{(3)}$} \\[2mm]
		\begin{tabular}{c|c|c|c}
			\hline
			$y$ & $u$ & $\left|\overline{G}^{(1)}-\overline{G}^{(3)}\right|$ & $DK2$ \\
			\hline
			0.10 & 0.10 &  0.0035 & 0.0779 \\
			& 0.25 &  0.0087 & \\
			& 0.50 &  0.0183 & \\
			& 1.00 &  0.0306 & \\
			& 2.00 &  0.0272 & \\
			&      &        & \\
			0.25 & 0.10 & 0.0094 & 0.0723 \\
			& 0.25 & 0.0156 & \\
			& 0.50 & 0.0244 & \\
			& 1.00 & 0.0322 & \\
			& 2.00 & 0.0252 & \\
			&      &        & \\
			0.50 & 0.10 & 0.0202 & 0.0592 \\
			& 0.25 & 0.0254 & \\
			& 0.50 & 0.0309 &  \\
			& 1.00 & 0.0324 & \\
			& 2.00 & 0.0217 & \\
			&      &        & \\
			1.00 & 0.10 & 0.0317 & 0.0413 \\
			& 0.25 & 0.0328 & \\
			& 0.50 & 0.0322 &  \\
			& 1.00 & 0.0270 &  \\
			& 2.00 & 0.0152 & \\
			&      &        & \\
			2.00 & 0.10 & 0.0227 & 0.0494 \\
			& 0.25 & 0.0211 & \\
			& 0.50 & 0.0183 &   \\
			& 1.00 & 0.0134 &  \\
			& 2.00 & 0.0068 & \\
			\hline
		\end{tabular}
	\end{minipage}
	\hfill
	\begin{minipage}{0.48\textwidth}
		\centering
		\scriptsize
		\textbf{(d) $\theta = 4$, $X^{(3)}$ vs $X^{(2)}$} \\[2mm]
		\begin{tabular}{c|c|c|c}
			\hline
			$y$ & $u$ & $\left|\overline{G}^{(3)}-\overline{G}^{(2)}\right|$ & $DK2$ \\
			\hline
			0.10 & 0.10 &  0.0001 & 0.0054 \\
			& 0.25 &  0.0003 & \\
			& 0.50 &  0.0007 & \\
			& 1.00 &  0.0016 & \\
			& 2.00 &  0.0023 & \\
			&      &        & \\
			0.25 & 0.10 & 0.0003 & 0.0052 \\
			& 0.25 & 0.0006 & \\
			& 0.50 & 0.0011 & \\
			& 1.00 & 0.0018 & \\
			& 2.00 & 0.0024 & \\
			&      &        & \\
			0.50 & 0.10 & 0.0008 & 0.0046 \\
			& 0.25 & 0.0011 & \\
			& 0.50 & 0.0015 &  \\
			& 1.00 & 0.0021 & \\
			& 2.00 & 0.0023 & \\
			&      &        & \\
			1.00 & 0.10 & 0.0017 & 0.0035 \\
			& 0.25 & 0.0019 & \\
			& 0.50 & 0.0021 &  \\
			& 1.00 & 0.0023 &  \\
			& 2.00 & 0.0020 & \\
			&      &        & \\
			2.00 & 0.10 & 0.0022 & 0.0027 \\
			& 0.25 & 0.0022 & \\
			& 0.50 & 0.0022 &   \\
			& 1.00 & 0.0020 &  \\
			& 2.00 & 0.0014 & \\
			\hline
		\end{tabular}
	\end{minipage}
	
\end{table}

\begin{rem}
	\label{remark-deficit}
	In Example $\ref{Example-deficit}$, the deviations observed between the bound DK2 and the exact value are natural since the bound is derived to hold uniformly with respect to $u$.
\end{rem}

\section{Ruin probability in a perturbed risk model}
\label{section-diffusion}
In the classical risk model the number of claims $N(t)$ in $\eqref{risk process}$ follows a homogeneous Poisson process with intensity $\lambda>0$. An extension of this is to add a diffusion term to account for additional uncertainties in the aggregate claims or the premium income. Gerber (1970) introduced the classical model perturbed by diffusion, with surplus process
\begin{equation}
	\label{surplus-diffusion}
	U(t)=u+ct-S(t)+\sigma B(t), \quad t\geq 0,
\end{equation}
where the dispersion parameter $\sigma>0$ and $\{B(t): t\geq 0\}$ is a standard Wiener process that is independent of the compound Poisson process $\{S(t):t \geq 0\}$ and of the individual claim sizes.\par
Dufresne and Gerber in \cite{DufresneGerber1991} studied three kinds of probabilities based on $\eqref{surplus-diffusion}:$ $\psi_d(u)=Pr(T<\infty, U(T)=0|U(0)=u)$ is the probability for ruin that is caused by oscillation, $\psi_s(u)=Pr(T<\infty, U(T)<0|U(0=u)$ is the probability that ruin is caused by a claim and $\psi_t(u)=Pr(T<\infty|U(0)=u)$, the probability of ruin.  We have that $\psi_t(u)=\psi_d(u)+\psi_s(u)$, with $\psi_d(0)=1$ and $\psi_s(0)=1$.\par
It was shown in \cite{DufresneGerber1991} showed that $\psi_t(u)=Pr(L^\ast>u)$ is the tail probability of the maximal aggregate loss $L^\ast=max\{u-U(t): t\geq 0\}$. The r.v. $L^\ast$ can be decomposed as 
\[
L^*=L_{o,0}+L_{c,1}+\cdots+L_{c,N}+L_{o,N}=\sum_{n=1}^{\infty}\left(L_{o,n-1}+L_{c,N}\right)+L_{o,N},
\]
with $L^*=L_{o,0}$ if $N=0$, where $L_o,N$ and $L_{c,N}$ are the amounts that result in the $(n+1)$-th and $n$-th record highs of the aggregate loss process $\{u-U(t)\}$ due to oscillation and a claim, respectively, and $N$ it the number of record highs of the process $\{u-U(t)\}$ caused by a claim. In addition, the r.v's $L_{o,0}, L_{o,1}, L_{o,2}...$ are identically distributed (as $L_o$) with common distribution function $H_1(u)=1-e^{-(c/D)u}$, where $D=\sigma^2/2$ and $L_{c,1}, L_{c,2}, L_{c,3}...$ are identically distributed (as $L_c$) with common distribution $H_2(u)=F_e(u)$. Also, $N, L_{o,0}, L_{c,1}, L_{o,1}, L_{c,2}, L_{o,2}...$ are independent. For further details and a probabilistic viewpoint of $\psi_t$ (see \cite{DufresneGerber1991} and \cite{Tsai2006}).\par
Furthermore, Tsai in \cite{Tsai2003} showed that
\begin{equation}
	\label{K(u)-equation}
	\overline{K}(u)=Pr(L_K>u)=\frac{1}{1+\theta}\psi_d(u)+\psi_s(u)=\sum_{n=1}^{\infty}\frac{\theta}{1+\theta} \left(\frac{\theta}{1+\theta}\right)^n \overline{A^{*n}}(u), \quad u\geq 0,
\end{equation}
with $\overline{K}(0)=1/(1+\theta)$, where $\overline{A}(x)=1-H_1 \ast H_2(x)=1-\int_{0}^{x}\overline{H_1}(x-t)dH_2(t)$ is the distribution function of $L_o+L_c$ with density $a(x)$. Since $L^\ast=L_K+L_o$, we have that $Pr(L^\ast>u)$, the ruin probability for surplus process $\eqref{surplus-diffusion}$, is a compound geometric convolution, given by
\begin{equation}
	\label{psit}
	\psi_t(u)=\overline{K\ast H_1}(u)=\sum_{n=0}^{\infty}\frac{\theta}{1+\theta} \left(\frac{\theta}{1+\theta}\right)^n \overline{A^{*n}\ast H_1}(u), \quad u\geq 0.
\end{equation}
Similarly, the expression for $\overline{K}(u)$ in $\eqref{K(u)-equation}$ can be considered as $\overline{K}(u)=Pr(L^K>u)$ where
\[
L^K=L_{o,0}+L_{c,1}+\cdots+L_{o,N-1}=\sum_{n=1}^{N}\left(L_{o, n-1}+L_{c,n}\right)
\]
where $L^K=0$ if $N=0$. When the diffusion term is removed (i.e., $\sigma=0$), then all $L_{o}$s disappear, implying model $\eqref{surplus-diffusion}$ reduce to the non-perturbed classical risk model and both $\psi_t(u)$ and $K(u)$ reduce to $\eqref{tail-G(x)}$. \par
\subsection{Continuity inequality for $K(u)$}
As seen in equation $\eqref{psit}$, $\psi_t(u)$ is the tail of the convolution of distribution function $K(u)$ and the diffusion term $H_1(u)$. Therefore, if an analytical expression for $\overline{K}(u)$ in $\eqref{K(u)-equation}$ is available, then $\psi_t(u)$ can also be obtained in explicit form. However, explicit expressions for $K(u)$ are generally unavailable, except in specific cases such as a combination of exponential claim distribution or a mixture of Erlangs (see \cite{Tsai2003}). Veraverbeke in \cite{Veraverbeke1993} investigated the asymptotic behavior of $\psi_t(u)$ in a classical risk model perturbed by diffusion and showed that the tail of $\psi_t(u)$ is related to the tail decay of $K(u)$ and $H_1(u)$. Moreover, an effective method to construct an upper bound for $\overline{K}(u)$ was proposed in \cite{Tsai2010}.  Hence, it is of interest to study the continuity properties of $K(u)$. In the following theorem, we derive a continuity inequality for the function in $\eqref{K(u)-equation}$.
\begin{thm}
	\label{K(u)-theorem}
	Let $\overline{K}, \overline{\widetilde{K}} \in \mathcal{K}_d$ be the tails of $\eqref{K(u)-equation}$ with $D\geq \widetilde{D}$ and $\mu \geq \widetilde{\mu}$. We also assume the net profit conditions $\phi, \widetilde{\phi}<1$. Then , we derive 
	\begin{align*}
		\sup_{u \geq 0}\left|\overline{K}(u)-\overline{\widetilde{K}} (u)\right|&\leq \frac{1}{c-\lambda \mu} \left[\lambda \mu\left(\frac{c}{D} \mathbb{K}(H_1, \widetilde{H}_1)+\frac{ |\widetilde{D}-D|}{D}+\frac{\mathbb{K}(F, \widetilde{F})}{\mu}  \right. \right. \\
		&\left. \left. +\frac{|\widetilde{\mu}-\mu|}{\mu}\right)+|\lambda\mu-\widetilde{\lambda}\widetilde{\mu}|\right]
	\end{align*}
\end{thm}

\textbf{Proof:}\\
	Let $\mathcal{K}_d$ be the space of all functions $x:[0, \infty)\rightarrow [0,1]$. The space $\mathcal{K}_d$ is a Banach space with the uniform metric $\nu_d(x, y):=\sup_{u \geq 0}|x(u)-y(u)|$, i.e. $(\mathcal{K}_d, \nu_d)$ is a complete metric space.\par
	It was proven in \cite{Tsai2003}  that $\overline{K}(u)$ satisfies the defective renewal equation:
	\begin{equation}
		\label{K(u)-equation-2}
		\overline{K}(u)=\frac{\lambda \mu}{c}\left(\overline{A}(u)+\int_{0}^{u}\overline{K}(u-t)dA(t)\right), \quad u \geq 0,
	\end{equation}
	For each $x\in \mathcal{K}_d$, we consider the operators $T_d, \widetilde{T}_d: \mathcal{K}_d \rightarrow \mathcal{K}_d$, defined by
	\begin{equation}
		\label{T-operator-Ku}
		T_d x(u)=\frac{\lambda \mu }{c} \left( \overline{A}(u)+\int_{0}^{u}x(u-t)dA(t)\right), \quad u\geq 0,
	\end{equation}
	\[
	\widetilde{T}_d x(u)=\frac{\widetilde{\lambda} \widetilde{\mu} }{c} \left( \overline{\widetilde{A}}(u)+\int_{0}^{u}x(u-t)d\widetilde{A}(t)\right), \quad u\geq 0.
	\]
	Given $u \geq 0$ and for all $x\in \mathcal{K}_d$ we have $T_d\mathcal{K}_d \subset \mathcal{K}_d$, since it follows immediately
	\begin{align*}
		T_d x(u)&= \frac{\lambda \mu}{c} \left(\overline{A}(u)+\int_{0}^{u}x(u-t)a(t)dt\right)\\
		&\leq \frac{\lambda \mu}{c} \left(\overline{A}(u)+\int_{0}^{u}a(t)dt\right) \leq \frac{\lambda \mu}{c}=\phi<1.
	\end{align*}
	Furthermore, these operators are contractive on $\mathcal{K}_d$ with modules $\phi$ and $\widetilde{\phi}$, respectively, since for all $u\geq 0$
	\begin{align*}
		\nu_d(T_d x, T_d y)	&=\frac{\lambda \mu}{c}\sup_{u \geq 0}\left|\int_{0}^{u}x(u-t)dA(t)-\int_{0}^{u}y(u-t)dA(t)\right|\\
		&\leq \frac{\lambda \mu}{c}\sup_{u \geq 0}\int_{0}^{u} \left|x(u-t)-y(u-t)\right|a(t)dt\\
		&\leq \frac{\lambda \mu}{c}\sup_{u \geq 0} \int_{0}^{u}\sup_{s\in [0,u]}\left|x(s)-y(s)\right| a(t)dt\\
		&\leq \frac{\lambda \mu}{c}\nu_d(x, y) \int_{0}^{\infty}a(t)dt=\frac{\lambda \mu}{c} \nu_d(x, y).
	\end{align*}
	Therefore, it follows that $\nu_d(T_d\overline{K}, T_d\overline{\widetilde{K}} )\leq \phi \nu_d(\overline{K}, \overline{\widetilde{K}} )$, where $\phi \leq 1$. Similarly, $\widetilde{T}_d$ is also a contractive operator on $\mathcal{K}_d$\\
	According to $\eqref{K(u)-equation-2}$, it is easy to show that $\overline{K}$ and $\widetilde{\overline{K}}$ are the unique fixed points of $\overline{K}=T_d\overline{K}$ and $\widetilde{\overline{K}}=\widetilde{T}_d\overline{\widetilde{K}} $. Hence, by triangle inequality it follows that
	\begin{align}
		\label{inequality-K1}
		v_{d}(\overline{K}, \overline{\widetilde{K}} )=v_{d}(T_d\overline{K}, \widetilde{T}_d\overline{\widetilde{K}} )& \leq v_{d}(T_d\overline{K}, T_d\overline{\widetilde{K}} )+v_{d}(T_d\overline{\widetilde{K}} , \widetilde{T}_d\overline{\widetilde{K}} ) \nonumber \\ 
		&\leq \frac{\lambda \mu}{c}  v_{d}(\overline{K}, \overline{\widetilde{K}} )+v_{d}(T_d \overline{\widetilde{K}}, \widetilde{T}_d \overline{\widetilde{K}} ) \nonumber \\ 
		&\leq \frac{c}{c-\lambda \mu} v_{d}(T_d \overline{\widetilde{K}}, \widetilde{T}_d\overline{\widetilde{K}}).
	\end{align}
	Furthermore,
	\begin{align}
		\label{inequality-K2}
		\nu_d(T_dx, \widetilde{T}_dx)
		&\leq \sup_{u \geq 0}\left| \frac{\lambda \mu}{c}\left( \overline{A}(u)+\int_{0}^{u}x(u-t) a(t)dt\right) \right. \nonumber \\
		& \left. - \frac{\lambda \mu}{c} \left(\overline{\widetilde{A}}(u)+\int_{0}^{u}x(u-t) \widetilde{a}(t)dt\right) \right| \nonumber \\
		& +\sup_{u \geq 0}\left|\frac{\lambda \mu - \widetilde{\lambda} \widetilde{\mu}}{c}\left(\overline{\widetilde{A}}(u)+\int_{0}^{u}x(u-t) \widetilde{a}(t)dt\right)\right| \nonumber \\
		& \leq \frac{\lambda \mu}{c}\int_{0}^{\infty}|a(t)-\widetilde{a}(t)|dt+\frac{|\lambda \mu - \widetilde{\lambda} \widetilde{\mu}|}{c} \int_{0}^{\infty}\widetilde{a}(t)dt,
	\end{align}
	where 
	\[
	a(t)=\int_{0}^{t}h_1(t-z) dH_2(z)=\int_{0}^{t}\frac{c}{D}e^{-(c/D)  (t-z)} \frac{\overline{F}(z)}{\mu}dz.
	\]
	It follows that  
	\begin{align*}
		\int_{0}^{\infty}|a(t)-\widetilde{a}(t)|dt	& \leq
		\int_{0}^{\infty}\left|\int_{0}^{t}\frac{c}{D}e^{-(c/D)(t-z)}\frac{\overline{F}(z)}{\mu}dz-\int_{0}^{t}\frac{c}{\widetilde{D}}e^{-(c/\widetilde{D})(t-z)}\frac{\overline{F}(z)}{\mu}dz \right| \\
		&+\left|\int_{0}^{t}\frac{c}{\widetilde{D}}e^{{-(c/\widetilde{D})(t-z)}}\frac{\overline{F}(z)}{\mu}dz -\int_{0}^{t}\frac{c}{\widetilde{D}}e^{-(c/\widetilde{D})(t-z)}\frac{\overline{\widetilde{F}}(z)}{\widetilde{\mu}}dz\right|dt  \\
		& \leq \int_{0}^{\infty} \frac{\overline{F}(z)}{\mu} \int_{z}^{\infty}\left|\frac{c}{D}e^{-(c/D)(t-z)}-\frac{c}{\widetilde{D}}e^{-(c/\widetilde{D})(t-z)}\right|dtdz \\
		&+ \int_{0}^{\infty} \left| \frac{\overline{F}(z)}{\mu}-\frac{\overline{\widetilde{F}}(z)}{\widetilde{\mu}} \right|\int_{z}^{\infty} \frac{c}{\widetilde{D}}e^{{-(c/\widetilde{D})(t-z)}}dtdz.
	\end{align*}
	Changing the integration variable by setting $t=y+z$, we obtain
	\begin{equation}
		\label{inequality-K3}
		\int_{0}^{\infty}|a(t)-\widetilde{a}(t)|dt \leq\frac{c}{D} \mathbb{K}(H_1, \widetilde{H}_1)+\frac{ |\widetilde{D}-D|}{D}+\frac{1}{\mu} \mathbb{K}(F, \widetilde{F})+\frac{|\widetilde{\mu}-\mu|}{\mu}.
	\end{equation}
	By inequalities $\eqref{inequality-K1}$, $\eqref{inequality-K2}$ and $\eqref{inequality-K3}$ we conclude the proof and the desired statement holds.
\hfil \qed \\

To illustrate the results of Theorem $\ref{K(u)-theorem}$ we obtain the following example:
\begin{exmp}
	For the surplus process $\eqref{surplus-diffusion}$, let  $X$ and $\widetilde{X}$ have the survival functions $\overline{F}_{X}(t)=e^{-3t}$ and $\overline{F}_{\widetilde{X}}(t)=(1/2)e^{-2t}+(1/2)e^{-6t}$, respectively, with $EX=E\widetilde{X}=1/3$, $\theta={\theta}=1,$ and $c={\widetilde{c}}=1$. For various choices of $D$ and $\widetilde{D}$, we can easily calculate the real values of the distance  $v_{d}(\overline{K}, \overline{\widetilde{K}})$, respectively (see Example 2 in \cite{Tsai2006}). Table $\ref{Ku-table-example}$ presents a comparison between the exact values of these distances and the results obtained from the bound $DK3$ in the Theorem $\ref{K(u)-theorem}$, with the last column showing the deviation of the ratio (Bound/Exact) from 1.
	
	\begin{table}[h!]
		\centering
		\begin{tabular}{ccccc}
			\hline
			$D$ & $\widetilde{D}$ & $\displaystyle \sup_{u \ge 0}|\overline{K}(u)-\overline{\widetilde{K}}(u)|$
			& Bound $DK3$ &  $|1-\text{Ratio}|$ \\ 
			\hline
			1 & 1/10 & 0.0854 & 0.4837 &  4.66 \\
			1/2 & 1/10 & 0.0559 & 0.4337 & 6.75  \\
			1/2 & 1/3 & 0.0271 & 0.2004 & 6.39  \\
			2 & 1 & 0.0496 & 0.2837 & 4.71  \\
			2 & 1/10 & 0.1148 & 0.5087 & 3.43  \\
			3 & 1/10 & 0.1305 & 0.5171 & 2.96  \\
			3 & 1/20 & 0.1334 & 0.5254 & 2.94  \\
			\hline
		\end{tabular}
		\caption{A mixture of two exponentials vs an exponential perturbed by diffusion.}
		\label{Ku-table-example}
	\end{table}
	
\end{exmp}

\subsection{Estimation for K(u)}
An approximation method for the ruin probability, $\psi(u)$, based on the contractive properties and the BFPT was proposed in \cite{SanchezBaltazar2020}. Following similar arguments, we use a suitable contracting operator $T_d$ on a certain Banach space $(\mathcal{K}_d, \nu_d)$, so that it can be applied to obtain an approximation of $\overline{K}$. \par
By the BFPT for contraction mappings, there exists a unique function (fixed point) $\overline{K} \in \mathcal{K}_d$, such that $T_d\overline{K}=\overline{K}$ and $\overline{K}$ is the limit of the functions 
\begin{equation}
	\label{equation--3}
	\overline{K}_n:=T_d\overline{K}_{n-1}=T_d^n \overline{K}_0, \quad n\geq 1,
\end{equation}
for some arbitrary function $\overline{K}_0$ of $\mathcal{K}_d$ and $T_d$ given in $\eqref{T-operator-Ku}$. Next, we consider the iterative sequence $\{\overline{K}_n\}_{n \geq 0}$ associated with $\overline{K}$ as defined in $\eqref{equation--3}$, 
\begin{equation}
	\label{equation--4}
	\overline{K}_n(u):= \frac{1}{1+\theta} \left(\overline{A}(u)+\int_{0}^{u}\overline{K}_{n-1}(u-x)dA(x)\right),
\end{equation}
where $\overline{K}_0 \in \mathcal{K}_d$. \par
Since $T_d$ is a contractive operator in Banach space $(\mathcal{K}_d, \nu_d)$, the BFPT (see \cite{GranasDugundji2003}) ensures that iterative sequence defined in $\eqref{equation--3}$ converges to the unique fixed point $\overline{K}$, that is
\[
\nu_d(\overline{K}_n, \overline{K})\rightarrow 0.
\]
The following result provides an explicit formula to compute the elements of the sequence $\{\overline{K}_n\}_{n \geq 0}$.
\begin{lm}
	\label{Lemma-2}
	Let $\overline{K}_n \in \mathcal{K}_d$, if $\overline{K}_0(u)=k \in (0,1)$ and for $u \geq 0$, then
	\[
	\overline{K}_n(u)= \left\{\begin{matrix}
		\phi - (1-k) \phi A(u),	&  n=1,\\
		\phi-(1-k) \phi^n A^{*(n)}(u)-(1-\phi) \sum_{i=1}^{n-1}\phi^i A^{*(i)}(u),	&  n\geq 2,\\
	\end{matrix}\right. 
	\]
	where $\phi=1/(1+\theta)<1$.
\end{lm}

\textbf{Proof:}\\
	The result can be demonstrated by mathematical induction. For $n=1$ applying $\eqref{equation--4}$ with $\overline{K}_0=k$, we have
	\begin{equation}
		\label{equation--5}
		\overline{K}_1(u)=TK_0(u)=Tk=\phi \left(\overline{A}(u)+\int_{0}^{u}k \cdot a(x)dx \right) =\phi - (1-k) \phi A(u).
	\end{equation}
	Now, suppose that the result holds for $n\geq 2$ with $u\geq 0$ fixed and for the next iteration we have
	\begin{align*}
		\overline{K}_{n+1}(u)&=T\overline{K}_n(u)\\
		&=\phi \left(\overline{A}(u)+ \int_{0}^{u}\left[\phi-(1-k)p^n A^{*n}(u-x) \right. \right. \\
		&\left. \left. -(1-\phi)\sum_{i=1}^{n-1}\phi^i A^{*i}(u-x)\right]a(x)dx\right)\\
		&=\phi \left(1-A(u)+\phi A(u)-(1-k) \phi^n A^{*(n+1)}(u) \right. \\
		&\left.-(1-\phi)\sum_{i}^{n-1}\phi^{i} A^{*(i+1)}(u)\right)\\
		&=\phi-(1-k)\phi^{n+1}A^{*(n+1)}(u) -(1-\phi) \sum_{j=1}^{n} \phi ^{j}A^{*(j+1)}(u).
	\end{align*}
	This completes the induction.  
\hfil \qed \\

\begin{rem}
	As $n\rightarrow \infty$, we reduce to equation $\eqref{K(u)-equation-2}$ as an immediate consequence of BFPT and Lemma $\ref{Lemma-2}$. In particular,  if $\overline{K}_0(u)=\phi$ and under the net profit conditions (see equation $\eqref{net profit}$) for classical model with $\overline{K}(u)$, then it follows that 
	\[
	\overline{K}(u)=\lim_{n \rightarrow \infty}\overline{K}_n(u)=(1-\phi)\sum_{i=1}^{\infty}\phi^i \overline{A^{*i}}(u), \quad u\geq 0,
	\]
	with $\overline{K}(0)=\phi$.
\end{rem}

\begin{exmp}
	Suppose that the individual claim amounts follow an exponential distribution with parameter $\beta$ in the surplus process in $\eqref{surplus-diffusion}$. We notice that $X$ satisfies the net profit condition if $\beta \in (\lambda/c, \infty)$. In this case the explicit solution for $\overline{K}(u)$ is given by 
	\begin{equation}
		\label{Ku-exact}
		K(u)=\theta \left(D_1 e^{-s_1 u}+D_2 e^{-s_2 u}\right)
	\end{equation}
	where $s_1$ and $s_2$ are the roots of the equation $s^2-(b_0+\beta)s+[\theta/(1+\theta)b_0 \beta=0$ with $b_0=c/D$. The corresponding constants are given by
	\[
	D_1=\frac{s_2}{\theta(1+\theta)\sqrt{(b_0-\beta)^2+4b_0 \beta/(1+\theta)}}>0,
	\]
	and 
	\[
	D_2=-\frac{s_1}{\theta(1+\theta)\sqrt{(b_0-\beta)^2+4b_0 \beta/(1+\theta)}}<0
	\]
	(see Example 2 in Tsai 2006).\par
	Therefore, we apply  BFPT and consider the corresponding iterative sequence of $\overline{K}(u)$, as defined in $\eqref{equation--4}$, in order to obtain an approximation of the function $\overline{K}(u)$. In the special case  $\beta=c/D$, where $A(u)$ follows an $Erlang(2, \beta)$, the following expression of $\overline{K}_n$ is obtained after some straightforward algebra and using mathematical induction,
	\begin{equation}
		\label{Kn-iterative}
		\overline{K}_n(u) =
		\begin{cases}
			\phi k + \phi (1-k) e^{-\beta u} (1+\beta u ), \quad n = 1, \\[0.3cm] 
			\phi e^{-\beta u} \mathcal{S}_{1}(\beta u)
			+ \sum_{m=2}^{n-1} \phi^{m} e^{-\beta u} \left[\mathcal{S}_{2m-1}(\beta u) - S_{2m-3}(\beta u)\right]
			\\
			+ \phi^{n}k \left(1 - e^{-\beta u}\mathcal{S}_{2n-3}(\beta u)\right)
			\\
			+ \phi^{n}(1-k) e^{-\beta u} \left[\mathcal{S}_{2n-1}(\beta u) - \mathcal{S}_{2n-3}(\beta u)\right], \quad n \geq 2,
		\end{cases}
	\end{equation}
	where $\mathcal{S}_m(z)=\sum_{r=0}^{m}\frac{z^r}{r!}$ (where $\mathcal{S}_{-1}(z)=0$ by convention).\par
	Tables \ref{tableKn} and \ref{tableKn-2} present the numerical approximations of $\overline{K}(u)$ obtained from the first five iterations of expression \eqref{Kn-iterative}, using the initial values $\overline{K}_0:=k=0.0, 0.1, 0.2, \ldots, 1.0$.
	Table \ref{tableKn} corresponds to the parameter set $\beta = 2$, $\lambda = c = 1/2$, $D = 1/4$, and $u = 1$, for which the exact value is $\overline{K}(1) = 0.3325717$.
	Table \ref{tableKn-2} refers to the case $\beta = 3/2$, $\lambda = 3/4$, $c = 2/3$, and $D = 4/9$, where the exact value equals $\overline{K}(1)=0.6573777$.
\end{exmp}

\begin{table}[h!]
	\centering
	{\normalsize 	\setlength{\tabcolsep}{4.5pt}
		\begin{tabular}{|c|c|c|c|c|c|c|}
			\hline
			$n$ & $k=0.0$  & $k=0.2$  & $k=0.4$  & $k=0.6$ &  $k=0.8$  & $k=1.0$ \\
			\hline
			$1$  & 0.2030029  & 0.2624023  & 0.3218018  & 0.3812012  & 0.4406006 &  0.5000000 \\
			$2$  & 0.3157823  & 0.3229262  & 0.3300700  & 0.3372138  & 0.3443576 &  0.3515015 \\
			$3$  & 0.3315714  & 0.3319855   & 0.3323996  & 0.3328137  & 0.3332278 &  0.3336419  \\
			$4$  & 0.3325381  & 0.3325518  & 0.3325655  & 0.3325793  & 0.3325930 &  0.3326067 \\
			$5$  & 0.3325709  & 0.3325712  & 0.3325717  & 0.3325720  & 0.3325721 &  0.3325724  \\
			\hline
		\end{tabular}
	}
	\caption{Approximations of the tail $\overline{K}(u)$ for the first five iterations when $X_i\sim Exp(2)$, $D=1, c=\lambda=1$ and $u=1$.}\
	\label{tableKn}
\end{table}

\begin{table}[h!]
	\centering
	{\normalsize 	\setlength{\tabcolsep}{4.5pt}
		\begin{tabular}{|c|c|c|c|c|c|c|}
			\hline
			$n$ & $k=0.0$  & $k=0.2$  & $k=0.4$  & $k=0.6$  & $k=0.8$  & $k=1.0$ \\
			\hline
			$1$  & 0.4183691  & 0.4846952  & 0.5510214  & 0.6173476  & 0.6836738 &  0.7500000 \\
			$2$  & 0.6301684  & 0.6375532  & 0.6449379  & 0.6523227  & 0.6597075 &  0.6670923 \\
			$3$  &  0.6559814  & 0.6563574  & 0.6567334  & 0.6571093  & 0.65745853 &  0.6578613 \\
			$4$  & 0.6573377  & 0.6573484  & 0.6573591  & 0.6573699  & 0.6573806 &  0.6573913 \\
			$5$  & 0.6573769  & 0.6573771  & 0.6573773  & 0.6573775  & 0.6573777 &  0.6573779 \\
			\hline 
		\end{tabular}
	}
	\caption{Approximations of the tail $\overline{K}(u)$ for the first five iterations when $X_i\sim Exp(3/2)$, $c=2/3, D=4/9, \lambda=3/4$ and $u=1$.}\
	\label{tableKn-2}
\end{table}

\section{Conclusions}
In this work, continuity inequalities were established for ruin-related quantities within the classical risk model. Specifically, Theorem $\ref{LKdistance}$ provides a continuity inequality for the ruin probability under claim distributions with finite $(\gamma+1)$-moments, while Theorem $\ref{theorem-deficit}$ presents a continuity estimate for the deficit at ruin. Both results were obtained using properties of contractive integral operators. \par
In Section $\ref{section-diffusion}$, the analysis was extended to the classical risk model perturbed by diffusion, where a continuity inequality for the supremum distance of the function $K(u)$ was derived. In this framework, although the ruin probability  $\psi_t(u)=\overline{K\ast H_1}(u)$, can be decomposed using the geometric r.v. $N$ to denote the number of recorded highs, the $\psi_t(u)$ itself is not a compound geometric distribution.  This observation motivates us to study the behavior of the compound geometric distribution $K(u)$. Consequently, we derived a continuity estimate for $K(u)$ and proposed an iterative approximation based on the Banach fixed point theorem.\par


\clearpage











\vspace{2cm}




\begin{thebibliography}{}
	
	\bibitem{AsmussenAlbrecher2010}
	{Asmussen, S., Albrecher, H.} (2010)
	\newblock {\em Ruin Probabilities, (2nd ed.)}.
	\newblock {World Scientific Publishing Co. Pte. Ltd., Singapore}.
	
	
	\bibitem{BarecheCherfaoui2019}
	{Bareche, A., Cherfaoui, M.} (2019)
	\newblock {Sensitivity of the Stability Bound for Ruin Probabilities
	to Claim Distributions}.
	\newblock {\em Methodology and Computing in Applied Probability}
	{\bf 21}, 1259-1281.


	\bibitem{BeirlantRachev1987}
	{Beirlant, J., Rachev, S.T.} (1987)
	\newblock {The problem of stability in insurance mathematics}.
	\newblock {\em Insurance: Mathematics and Economics} {\bf 6}, 179-188.
	
	\bibitem{BenouaretAissani2010}
	{Benouaret, Z., A\"{i}ssani, D.} (2010)
	\newblock {Strong stability in a two-dimensional classical risk model with independent claims}.
	\newblock {\em Scandinavian Actuarial Journal} {\bf 2}, 83-92.
	

	
	\bibitem{DufresneGerber1991}
	{Dufresne, F. and Gerber, H.U.} (1991).
	\newblock {Risk theory for the compound Poisson process that is perturbed by diffusion}.
	\newblock {\em Insurance: Mathematics and Economics} 10:51-59.


	\bibitem{EnikeevaKalashnikov2001}
	{Enikeeva, F., Kalashnikov, V., Rusaityte} (2001)
	\newblock {Continuity estimates for ruin probabilities}.
	\newblock {\em Scandinavian Actuarial Journal} {\bf 1}, 18-39.

	\bibitem{GajekRudz2013}
	{Gajek, L., Rud\'z, M.} (2013)
	\newblock {Sharp approximations of ruin probabilities in the discrete time models}.
	\newblock {\em Scandinavian Actuarial Journal} {\bf 5}, 352-382.

	\bibitem{GajekRudz2018}
	{Gajek, L., Rud\'z, M.} (2018)
	\newblock {Banach contraction principle and ruin probabilities in regime-switch models}.
	\newblock {\em Insurance: Mathematics and Economics} {\bf 80}, 45-53.
	
	\bibitem{GajekRudz2025}
	{Gajek, L., Rud\'z, M.} (2025)
	\newblock {Applications of the Banach fixed-point theorem to analyze insolvency problems of an insurance company}.
	\newblock {\em  Journal of Fixed Point Theory and Applications}, {\bf 27}, 34.

	\bibitem{Gerber1970}
	{Gerber, H.U.} (1970).
	\newblock {An extension of the renewal equation and its application in the collective theory of risk}.
	\newblock {\em Skandinavisk Aktuarietidskrift}, 205-210.


	\bibitem{GerberGoovaertsKaas1987}
	{Gerber, H.U., Goovaerts, M.J., Kaas, R.} (1987).
	\newblock {On the probability and severity of ruin}.
	\newblock {\em ASTIN Bulletin} 39:151-163.


	\bibitem{GordienkoOrtega2016}
	{Gordienko, E., V\'azquez-Ortega, P.} (2016)
	\newblock {Simple continuity inequalities for ruin probability in the classical risk model}.
	\newblock {\em ASTIN Bulletin} {\bf 46}, 801--814.
	
	\bibitem{GordienkoOrtega2018}
	{Gordienko, E., V\'azquez-Ortega, P.} (2018)
	\newblock {Continuity inequalities for multidimensional renewal risk models}.
	\newblock {\em Insurance: Mathematics and Economics} {\bf 82}, 48--54.
	
	
	\bibitem{GranasDugundji2003}
	{Granas, A., Dugundji} (2003)
	\newblock {\em Fixed Point Theory}.
	\newblock {Springer-Verlag, New York}.
	
	
	\bibitem{GrubelPitts1993}
	{Gr\"ubel, R, Pitts, S. M.} (1993)
	\newblock {Nonparametric estimation in renewal theory I: the empirical renewal function}.
	\newblock {\em Annals Statistics}  {\bf 21}, 1431-1451.

	\bibitem{HarfoucheBareche2021}
	{Harfouche, Z., Bareche, A.} (2021)
	\newblock {Semi-parametric approach for approximating the ruin probability of classical risk models with large claims}.
	\newblock {\em Communications in Statistics - Simulation and Computation}  {\bf 11}, 5585-5604.

	\bibitem{Jiang2021}
	{Jiang Z.} (2021)
	\newblock {Banach contraction principle, q-scale function and ultimate ruin probability under a Markov-modulated classical risk model}.
	\newblock {\em Scandinavian Actuarial Journal} {\bf 3}, 234-243.
	
	
	\bibitem{Kalashnikov1997}
	{Kalashnikov, V.V.} (1997)
	\newblock {\em Geometric sums: Bounds for Rare Events With Applications: Risk Analysis, Reliability, Queuing}.
	\newblock {Kluwer, New York}.
	
	\bibitem{Kalashnikov2000}
	{Kalashnikov, V.V.} (2000)
	\newblock {The stability concept for stochastic risk models}.
	\newblock {\em Working paper No.166. Laboratory of Actuarial Mathematics, University of Copenhagen}.
	
	\bibitem{PittsPolitis2008}
	{Pitts, M., Politis, K.} (2008)
	\newblock { Approximations for the moments of ruin time in the compound poisson model }.
	\newblock {\em Insurance: Mathematics and Economics 42, 668-679}.
	
	\bibitem{Politis2006}
	{Politis, K.} (2006)
	\newblock { A functional approach for ruin probabilities}.
	\newblock {\em Stochastic Models 22, 509-536}.
	
	
	\bibitem{Pollard1984}
	{Pollard, D.} (1984)
	\newblock {\em Convergence of Stochastic Processes}.
	\newblock {Springer, New York}.
	


	
	
	
	
	\bibitem{Rusaityte2001}
	{Rusaityte, D.} (2001)
	\newblock {Stability bounds for ruin probabilities in a Markov modulated risk model with investments}.
	\newblock {\em Working Paper Nr 178. Lab. of Actuarial
		Mathematics, University of Copenhagen.}


	
	\bibitem{SanchezBaltazar2020}
	{S\'anchez, J.M., Baltazar-Larios, F.} (2020)
	\newblock {Approximations of the ultimate ruin probability in the classical risk model using the Banach's Fixed-Point Theorem and continuity of the ruin probability}.
	\newblock {\em Kybernetika}
	{\bf 58}(2), 254-276.
	

	\bibitem{SantanRincon2020}
	{Santana, D.J., Rinc\'on} (2020)
	\newblock {Approximations of the ruin probability in a discrete time risk model}.
	\newblock {\em Modern Stochastics: Theory and Applications}
	\{7\}(3), 221-243.
	
	
	\bibitem{TouaziBenouaret2017}
	{Touazi, A., Benouaret, Z., A\"issani, D., Adjabi, S.} (2017)
	\newblock {Non-parametric estimation of the claim amount in the strong stability analysis of the classical risk model}.
	\newblock {\em Insurance: Mathematics and Economics} 74, 78-83.
	
	
	\bibitem{Tsai2003}
	{Tsai, C.C.L.} (2003)
	\newblock {On the expectations of the present values of the time of ruin perturbed by diffusion}.
	\newblock {\em Insurance: Mathematics and Economics} \textbf{32}(3), 413-429.
	
	\bibitem{Tsai2006}
	{Tsai, C.C.L.} (2006)
	\newblock {On the stop-loss transform and order for the surplus process perturbed by diffusion}.
	\newblock {\em Insurance: Mathematics and Economics} \textbf{39}(1), 151-170.
	
	\bibitem{Tsai2010}
	{Tsai, C.C.L.} (2010)
	\newblock {An effective method for constructing bounds for ruin probabilities for the surplus process perturbed by diffusion}.
	\newblock {\em Scandinavian Actuarial Journal} 2010, 3, 200-220.
	
	
	\bibitem{Veraverbeke1993}
	{Veraverbeke, N.} (1993)
	\newblock {Asymptotic Estimations for the Probability of Ruin in a Poisson Model with diffusion}.
	\newblock {\em Insurance: Mathematics and Economics} {\bf 13}, 57-62.
	
	
	\bibitem{Willmot2002}
	{Willmot, G. E.} (2002)
	\newblock {Compound geometric residual lifetime distributions and the deficit at ruin}.
	\newblock {\em Insurance: Mathematics and Economics} 30(3), 421-438.
	
	\bibitem{WillmotLin1998}
	{Willmot, G. E., Lin, X. S.} (1998)
	\newblock {Exact and approximate properties of the distribution of surplus before and after ruin}.
	\newblock {\em Insurance: Mathematics and Economics} 23(1), 91-110.
	
	
	

	
\end{thebibliography}
\end{document}